\documentclass[11pt,letterpaper,reqno]{amsart}

\usepackage[margin=1.02in]{geometry}
\usepackage[T1]{fontenc}
\usepackage[utf8]{inputenc}
\usepackage{lmodern}
\usepackage{microtype}
\usepackage{tikz-cd}
\usepackage{verbatim}

\usepackage{amssymb,mathtools,mathrsfs}
\usepackage{graphicx,booktabs,enumitem}

\usepackage{xcolor}
\definecolor{paperblue}{HTML}{124878}
\usepackage[
  colorlinks=true,
  linkcolor=paperblue,
  citecolor=paperblue,
  urlcolor=paperblue,
  linktoc=all,
  bookmarksnumbered=true,
  pdfencoding=auto,
  psdextra
]{hyperref}
\usepackage[
  backend=biber,
  style=alphabetic,
  maxnames=10,
  giveninits=true,
  sortcites=true,
  sorting=nyt          
]{biblatex}

\AtEveryBibitem{%
  \ifboolexpr{
    test {\iffieldundef{year}}
    and
    test {\iffieldundef{issue}}
  }
    {\renewbibmacro*{issue+date}{}}
    {}%
}

\usepackage{aliascnt}
\usepackage[capitalize,nameinlink,noabbrev]{cleveref}

\newcommand{\C}{\mathbb{C}}

\newcommand{\R}{\mathbb{R}}

\newcommand{\Z}{\mathbb{Z}}
\newcommand{\N}{\mathbb{N}}
\newcommand{\Pb}{\mathbb{P}}

\DeclareMathOperator{\supp}{supp}

\DeclareMathOperator{\image}{Im}
\DeclareMathOperator{\codim}{codim}
\DeclareMathOperator{\vol}{vol}

\renewcommand{\d}{\mathop{}\!\mathrm{d}}
\newcommand{\ddc}{\mathop{}\!\mathrm{d}\mathrm{d}^c}

\newcommand{\cB}{\mathcal{B}}
\newcommand{\cC}{\mathcal{C}}

\newcommand{\cE}{\mathcal{E}}

\newcommand{\cI}{\mathcal{I}}

\newcommand{\cO}{\mathcal{O}}

\newcommand{\lp}{\langle}
\newcommand{\rp}{\rangle}

\newcommand{\id}{\mathrm{id}}
\newcommand{\FuStu}{\omega_{\rm{FS}}}
\newcommand{\tildeGamma}{\widetilde{\Gamma}}

\newcommand{\into}{\hookrightarrow}

\newcommand{\newsharedtheorem}[2]{%
  \newaliascnt{#1}{theorem}%
  \newtheorem{#1}[#1]{#2}%
  \aliascntresetthe{#1}%
}

\theoremstyle{plain}
\newtheorem{theorem}{Theorem}[section]
\newsharedtheorem{lemma}{Lemma}
\newsharedtheorem{proposition}{Proposition}
\newsharedtheorem{corollary}{Corollary}

\theoremstyle{definition}
\newsharedtheorem{definition}{Definition}
\newsharedtheorem{example}{Example}

\theoremstyle{remark}
\newsharedtheorem{remark}{Remark}

\Crefname{theorem}{Theorem}{Theorems}
\Crefname{lemma}{Lemma}{Lemmas}
\Crefname{proposition}{Proposition}{Propositions}
\Crefname{corollary}{Corollary}{Corollaries}
\Crefname{definition}{Definition}{Definitions}
\Crefname{example}{Example}{Examples}
\Crefname{remark}{Remark}{Remarks}
\Crefname{appendix}{Appendix}{Appendices}

\numberwithin{equation}{section}

\title[Geometric structure of Green currents for holomorphic correspondences]{Geometric structure of Green currents for holomorphic correspondences on projective surfaces}
\author{Muhan Luo}
\address{Department of Mathematics, National University of Singapore, Singapore 119076}
\curraddr{Université d'Orléans, Université de Tours, CNRS, IDP, UMR 7013, 45067 Orléans, France}
\email{luomuh98@gmail.com}

\author{Marco Vergamini}
\address{Scuola Normale Superiore, Piazza dei Cavalieri 7, 56126 Pisa, Italy}
\email{marco.vergamini@sns.it}

\date{\today} 
\subjclass[2020]{37F80, 32U40, 32H50}
\keywords{Holomorphic correspondences, Green currents, woven currents, dynamical degrees}

\begin{document}

\begin{abstract}
 Let $f$ be a holomorphic correspondence on a complex projective surface $X$ such that $f^{-1}$ is also holomorphic. Assume that the action of $f^*$ on cohomology is simple and that the dynamical degrees of $f$ satisfy $d_1>\max\{d_0,d_2\}$. We prove that the Green currents $T^\pm$ associated with $f$ are woven. We also provide a large family of examples of holomorphic correspondences satisfying our assumptions, obtained as generic complete intersections in $X\times X$.
\end{abstract}

\maketitle
\tableofcontents

\section{Introduction and main result}

The geometric structure of invariant currents is an important aspect of higher-dimensional holomorphic dynamics. Bedford, Lyubich, and Smillie introduced \textit{laminar currents} in their study of polynomial automorphisms of $\mathbb C^2$ \cite{BLS93a}. Roughly speaking, a positive $(1,1)$-current on a complex surface is laminar if it is locally obtained by averaging integration currents over compatible holomorphic curves: the curves may intersect only along subsets that are relatively open in both. A weaker geometric notion is that of a \textit{woven current}, for which no such compatibility condition is imposed. We refer to \Cref{woven_and_criterion} for the definition used in this paper.

The notions of laminar and woven currents have been used to study intersections of currents and the dynamical properties of invariant measures; see for example \cite{Dujardin04PubM,DDG3ENS,Dujardin06duke}. Laminar or woven structures of Green currents have been established for meromorphic maps of projective surfaces \cite{DDG1Indiana,Dujardin03MAnn}, automorphisms of projective K3 surfaces \cite{Cantat01Acta} and compact K\"ahler manifolds \cite{dTH-Dinh12AM,D05JGA}, and H\'enon-like maps \cite{BDR26JIMJ,LZ25}. See also \cite[Section~1]{Dujardin22} for an overview of geometric currents and their applications in holomorphic dynamics.

In this paper we address similar questions in the setting of holomorphic correspondences, which roughly speaking can be considered as multivalued holomorphic functions. There have been lots of recent works on the dynamics of correspondences. Their Green currents and related equidistribution results have been studied in \cite{DNV18AM,AV24AMP,LV26}. In particular, under the existence of a gap between two consecutive dynamical degrees, the authors constructed Green currents and established regularity properties of their super-potentials \cite{LV26}. On the other hand, the multivalued nature of correspondences causes significant difficulties when studying the geometric structure of the Green currents, which has not been addressed to the best of our knowledge. In this paper, we address this question for correspondences on projective surfaces.

\medskip

Let $X$ be a smooth complex projective surface, and let $f$ be a holomorphic correspondence on $X$ such that $f^{-1}$ is also holomorphic. Thus the \textit{graph} of $f$, denoted by $\Gamma$, is an effective $2$-cycle in $X\times X$ whose irreducible components project surjectively and with finite fibres onto both factors. For every $0\leq q\leq2$, we denote by $d_q$ the \textit{dynamical degrees of order $q$} of $f$, which is the spectral radius of $f^*$ acting on $H^{q,q}(X,\mathbb R)$. In particular, $d_0$ and $d_2$ are the generic numbers of images and inverse images of a point, respectively, counted with multiplicity. We consider the case
\[
    d_1>\max\{d_0,d_2\}.
\]
We also assume that the action of $f^*$ on $H^{1,1}(X,\mathbb R)$ is \textit{simple}, meaning that $d_1$ is a simple eigenvalue and the unique eigenvalue of maximal modulus.

By \cite[Proposition~4.4]{DS07MM}, the pullback and pushforward of $f$ on positive closed $(1,1)$-currents are well-defined in this setting. By \cite[Theorem~1.1]{LV26}, there are positive closed $(1,1)$-currents $T^\pm$, called the \textit{Green currents} of $f$, satisfying
\[
    f^*T^+=d_1T^+\qquad \text{and }\qquad f_*T^-=d_1T^-.
\]
Our main result is the following.

\begin{theorem}\label{thm:main}
Let $f$ be a holomorphic correspondence on a smooth complex projective surface $X$, and assume that $f^{-1}$ is also holomorphic. Suppose that the action of $f^*$ on $H^{1,1}(X,\mathbb R)$ is simple and that $d_1>\max\{d_0,d_2\}$. Then the Green currents $T^\pm$ associated with $f$ are woven.
\end{theorem}

It suffices to prove the result for $T^+$, since the statement for $T^-$ follows by applying the same argument to $f^{-1}$. Indeed, the dynamical degrees $d_q'$ of $f^{-1}$ satisfy $d_q'=d_{2-q}$, and $(f^{-1})^*=f_*$ is adjoint to $f^*$ on $H^{1,1}(X,\mathbb R)$, so the simplicity assumption is preserved. The strategy of our proof can be summarized as follows. We follow the main idea of \cite{Dujardin03MAnn,DDG1Indiana}. However, since $f$ is multivalued, some essential parts of the original proof have to be replaced. 

Fix an embedding $X\hookrightarrow\mathbb P^N$, and write $H_X:=H\cap X$ for a hyperplane $H\subset\mathbb P^N$. We first show that, for almost every $H$, the pullback cycles $C_n:=(f^n)^*(H_X)$ and their currents of integration $[C_n]$ satisfy
\[
    d_1^{-n}[C_n]\longrightarrow cT^+,\qquad\deg C_n\asymp d_1^n,
\]
where $c>0$ is independent of $H$. This is done by adapting the argument of \cite{DDG1Indiana}. By the geometric criterion in \Cref{genus_to_woven}, it remains to prove the bound $g(C_n)=O(d_1^n)$, where the geometric genus $g(C_n)$ is computed using the normalizations of the irreducible components of $C_n$, counted with multiplicity.

We proceed to construct smooth, possibly disconnected surfaces $\widetilde\Gamma_n$. These surfaces are equipped with maps $\Pi_{1,n},\Pi_{2,n}:\widetilde\Gamma_n\to X$ and satisfy
\[
    (\nu_n)_*\widetilde\Gamma_n=\Gamma_n, \quad\text{for }\nu_n:=(\Pi_{1,n},\Pi_{2,n}),
\]
where $\Gamma_n$ is the graph cycle of $f^n$, along with a commutative diagram \eqref{large-pyramid}. An important distinction with respect to the single-valued setting is that the maps $\Pi_{1,n}$ and $\Pi_{2,n}$ need not be generically one-to-one. While the idea of \cite{Dujardin03MAnn,DDG1Indiana} still applies, its application here requires a non-trivial adaptation to account for these technical changes.

For a suitable generic hyperplane $H$, the curves $\widehat C_n:=\Pi_{2,n}^*(H_X)$ are smooth and possibly disconnected, and satisfy
\[
    C_n=(\Pi_{1,n})_*\widehat C_n\qquad \text{and} \qquad \widehat C_n=\beta_n^*\widehat C_{n-1},
\]
where the map $\beta_n$ is constructed in \Cref{subsec:construct-model}. Suppose we have the decomposition $\widehat{C}_n=\sum_j \widehat{C}_{n,j}$ where each $\widehat{C}_{n,j}$ is an irreducible component. Denote by $G_n:=\sum_j g(\widehat C_{n,j})$. Applying the Riemann--Hurwitz formula componentwise gives
\[
    G_n\leq d_2G_{n-1}+O(d_1^{n-1}).
\]
From this and the hypothesis $d_2<d_1$, we can deduce that $G_n=O(d_1^n)$. Note that an irreducible component $\widehat C_{n,j}$ does not necessarily need to be a normalization of its image $C_{n,j}=\Pi_{1,n}(\widehat C_{n,j})$. Thanks to the universal property of normalizations, for each irreducible component $\widehat C_{n,j}$, the restriction of $\Pi_{1,n}$ factors through a finite map
\[
    \phi_{n,j}:\widehat C_{n,j}\longrightarrow C_{n,j}^{\nu},
\]
where $C_{n,j}^{\nu}$ is the normalization of $C_{n,j}$. The degree of $\phi_{n,j}$ may be greater than one. A second application of the Riemann--Hurwitz formula, to the maps $\phi_{n,j}$, gives the desired estimate $g(C_n)=O(d_1^n)$. 

\begin{remark}
    Our argument is developed in dimension two, making use of the genera of curves and the Riemann--Hurwitz formula. In higher dimension, geometric structures of currents have been studied by lifting to Grassmannian bundles; see \cite{BDR26JIMJ,dTh,dTH-Dinh12AM,DS16MAnn,LZ25}. It would be interesting to adapt such methods to holomorphic correspondences, taking into account the singularities of the iterated fibre products and the multiplicities arising from composition.
\end{remark}

\noindent\textbf{Outline of the paper.} The paper is organized as follows. \Cref{prels} recalls the necessary background and proves the approximation by pullbacks of hyperplanes. \Cref{iterated_graphs} constructs the compatible smooth models and the approximating curves. The proof of \Cref{thm:main} is completed in \Cref{proof_section}. Finally, in \hyperref[sec:example]{Appendix}, we give a family of examples as generic complete intersections in $X\times X$.

\medskip

\noindent\textbf{Acknowledgements.} This project was mainly carried out during a visit of the second author to the Laboratoire de Probabilités, Statistiques et Modélisation at the Sorbonne Université in Paris. He would like to thank the Laboratoire and the university for their hospitality, and in particular Professor Romain Dujardin for many useful discussions. He would also like to thank Professor Davide Lombardo for useful discussion concerning some algebraic geometry aspects of \Cref{prels,iterated_graphs}. The first author has received funding from the National University of Singapore and MOE of Singapore through the grant A-8002488-00-00 and from the ANR grant DynAtrois.

\medskip

\noindent\textbf{Declaration of AI use.} The mathematical content of the main part of the paper is due to the authors. GPT-5.6 was used to check and refine the proofs in \cref{iterated_graphs,proof_section}. The examples in \hyperref[sec:example]{Appendix} were found with the assistance of GPT-5.6. The authors take full responsibility for all mathematical content.

\section{Preliminaries} \label{prels}

\subsection{Cycles and intersection theory}

We recall some facts about cycles and intersection products that will be used below; see \cite[Chapter 0, Section 4 and Chapter 1, Section 1]{GH94PoAG}. For simplicity, we restrict our discussions to a complex projective manifold $M$.

A \emph{$k$-cycle} on $M$ is a finite sum $C=\sum_j m_jV_j$, where the $V_j$'s are irreducible closed subvarieties of dimension $k$ and $m_j\in\Z$. The cycle is \emph{effective} if all its coefficients are nonnegative. We write $C\ge C'$ if $C-C'$ is effective. A \emph{divisor} is a cycle of codimension one. Throughout the paper, we identify a pure-dimensional analytic set with the associated cycle it defines. For a cycle $C$, we denote its support by $|C|$ and its current of integration by $[C]$.

Two irreducible subvarieties $V,W\subset M$ \emph{intersect properly} if every irreducible component of $V\cap W$ has dimension $\dim V+\dim W-\dim M$; the empty intersection is also allowed. In this case, their \emph{intersection product} is the cycle
\[
    V\cdot W=\sum_Z i_Z Z,
\]
where $Z$ ranges over the irreducible components of $V\cap W$ and $i_Z$ denotes the intersection multiplicity. This can be extended bi-linearly to cycles whose components intersect properly. 

Let $h:Y\to X$ be a proper holomorphic map between complex manifolds, and let $B$ be a cycle of pure codimension $r$ on $X$. Assume that, for every irreducible component $W$ of $B$, the inverse image $h^{-1}(W)$ is empty or has pure codimension $r$ in $Y$, so that the pullback $h^*B$ is well defined. The \textit{projection formula} gives
\begin{equation}\label{eq:cycle-projection}
    h_*(h^*B)=B\cdot h_*Y
\end{equation}
as an equality of cycles.

When an intersection cycle is zero-dimensional, we also denote by $C\cdot C'$ its \textit{numerical intersection product}, i.e., the sum of its coefficients. In particular, let $h:Y\to X$ be a dominant generically finite holomorphic map between smooth connected compact complex surfaces. For divisors $A,B$ on $X$ and $D$ on $Y$, the numerical intersection products satisfy
\[
    h^*A\cdot h^*B=d_t(h)(A\cdot B), \qquad h_*D\cdot A=D\cdot h^*A,
\]
where $d_t(h)$ is the topological degree of $h$. 

We denote by $R_h$ the \emph{ramification divisor}, i.e., the divisor locally defined by the holomorphic function $\text{Jac}(h)$, the Jacobian determinant. We say that $h$ is \emph{unramified} if $\text{Jac}(h)$ never vanishes. By the holomorphic inverse function theorem, this is equivalent to saying that $h$ is a local biholomorphism.

\subsection{Fibre product and base change.} We recall the notions of fibre product and base change, as well as some preservation properties that will be used below.

\begin{definition}
    Let $g:X\rightarrow Z$ and $h:Y\rightarrow Z$ be two holomorphic functions between complex spaces. The \emph{fibre product} of $g$ and $h$ is the analytic subset of $X\times Y$ given by
    $$X\times_Z Y:=\{(x,y) \in X\times Y\mid g(x)=h(y)\}.$$
    Let $P_1$ and $P_2$ the projection of $X\times Y$ on $X$ and $Y$, respectively. Then $P_2|_{X\times_Z Y}\rightarrow Y$ is the \emph{base change} of $g$ and $ P_1|_{X\times_Z Y}\rightarrow X$ is the \emph{base change} of $h$.
\end{definition}

\begin{lemma} \label{base-change}
    Let $g:X\rightarrow Z$ and $h:Y\rightarrow Z$ be two holomorphic functions between complex spaces, and let $X\times_Z Y$ be their fibre product.
    \begin{enumerate}
        \item If $g$ is finite, then its base change is finite of the same degree.
        \item If $g$ is unramified, then its base change is unramified.
    \end{enumerate}

    In particular, if $g$ is a finite unramified covering, then $P_2|_{X\times_Z Y}$ is a finite unramified covering of the same degree.
\end{lemma}

\begin{proof}
    Suppose that $g$ is finite. Take a point $y\in Y$. By definition of fibre product, we have
    $$P_2|_{X\times_Z Y}^{-1}(y)=\{(x,y)\mid g(x)=h(y)\}.$$
    Counting multiplicities, this set has the same cardinality of $g^{-1}(h(y))$. This proves that the base change of $g$ has finite fibres, and its degree is the same of $g$. We have to prove that it is proper. Consider a compact set $K \subseteq Y$. We have
    $$P_2|_{X\times_Z Y}^{-1}(K)=\{(x,y) \in X\times Y\mid g(x)=h(y), y\in K\}=\big(g^{-1}(h(K))\times K\big)\cap (X\times_Z Y).$$
    Since $h$ is continuous, $h(K)$ is compact. Since $g$ is a finite holomorphic map, it is proper, so $g^{-1}(h(K))$ is also compact and $g^{-1}(h(K))\times K$ is compact too. Its intersection with $X\times_Z Y$, which is analytic hence closed, is a closed subset of $g^{-1}(h(K))\times K$ itself. Therefore, it is compact. This implies that the base change of $g$ is proper, and since it has finite fibres it is finite. This proves point (1).

    Suppose now that $g$ is unramified. Take $(x,y)\in X\times_Z Y$, and consider neighbourhoods $U$ of $x$ and $V$ of $g(x)$ such that $g|_U:U\rightarrow V$ is a homeomorphism. Consider then the open sets $V'=h^{-1}(V)$ and $U':=P_2|_{X\times_Z Y}^{-1}(V')\cap(U\times Y)=(X\times_ZY)\cap(U\times V')$. The base change $P_2|_{X\times_Z Y}$ restricted to $U'$ is a homeomorphism between $U'$ and $V'$ with inverse $(g|_U^{-1}\circ h,\id)$. Therefore, the base change is unramified. This proves point (2).

    Finally, if $g$ is a finite unramified covering of degree $d$, the first two assertions show that $P_2|_{X\times_Z Y}$ is finite and unramified. Moreover, for every $y\in Y$, the fibre $P_2|_{X\times_Z Y}^{-1}(y)$ is naturally identified with $g^{-1}(h(y))$, which consists of $d$ points. Hence, $P_2|_{X\times_Z Y}$ has degree $d$. This concludes the proof
\end{proof}

\subsection{Holomorphic correspondences and dynamical degrees}\label{subsec:correspondence}
We recall some basic definitions and operations for holomorphic correspondences; see \cite[Section~3]{DS06CMH} for the general framework of meromorphic transformations.

Let $X$ be a smooth complex projective surface, i.e., such that there exists an embedding $i: X\into \Pb^N$ for some $N$ large enough. From now on, we fix such an $i$. The K\"ahler form on $X$ can be chosen to be $\omega:=i^*(\FuStu)$. For a generic hyperplane $H$ in $\Pb^N$, denote by $H_X:=H\cap X$ the hyperplane section in $X$. It defines a positive closed $(1,1)$-current $[H_X]$ on $X$.

Let $\pi_1$ and $\pi_2$ be the canonical projection maps from $X\times X$ to its two factors. Let $\Gamma=\sum_j\gamma_j$ be an effective two-cycle, that is, a finite sum of irreducible analytic sets $\gamma_j$ in $X\times X$. We assume that $\pi_1(\gamma_j)=\pi_2(\gamma_j)=X$ for all $j$. The {\it meromorphic correspondence} $f$ associated to $\Gamma$ is given by 
\[
    f(A):=\pi_2(\pi_1^{-1}(A)\cap\Gamma)
\]
for any set $A\subset X$. The cycle $\Gamma$ is called the {\it graph} of $f$. We say that $f$ is a {\it holomorphic correspondence} if the {\it indeterminacy set} $I(f):=\{x\in X: \dim\pi_1^{-1}(x)\cap \Gamma>0\}$ is empty.

Given two holomorphic correspondences $f$ and $g$, the graph of their composition can be described as follows. Suppose first that their graphs $\Gamma_f$ and $\Gamma_g$ are irreducible subvarieties. The product $\Gamma_g\times\Gamma_f$ is naturally included in $X^4= \{(x_1,x_2,x_3,x_4):x_i \in X\}$. Set
\[
    D:=X\times\Delta_X\times X=\left\{(x_1,x_2,x_3,x_4)\in X^4:x_2=x_3\right\},
\] 
Define $\widehat{\Gamma}_{f\circ g}:=(\Gamma_g \times\Gamma_f)\cap D$, and let $\Gamma_{f\circ g}$ be the push forward of $\widehat{\Gamma}_{f\circ g}$ in the sense of cycles by the projection
\[
    p_{1,4}: (x_1,x_2,x_3,x_4)\to (x_1,x_4).
\]
Then $f\circ g$ is by definition the correspondence whose graph is $\Gamma_{f\circ g}$. In general, we write $\Gamma_f =\sum_i \Gamma_f^i$ and $\Gamma_g =\sum_j \Gamma_g^j$, where $\Gamma_f^i$ and $\Gamma_g^j$ are irreducible, and extend the above definition by linearity. The composition is associative and in particular, we can consider iterates of $f$. 

More precisely, we have the following, which will be used in \Cref{iterated_graphs}.

\begin{lemma}\label{lem:composition-cycle}
Let $f$ and $g$ be holomorphic correspondences on $X$. Denote their graphs by $\Gamma_f$ and $\Gamma_g$, respectively. Then the intersection cycle
\[
    \mathcal C(f,g):=D\cdot(\Gamma_g\times\Gamma_f)
\]
is well-defined and has pure dimension $2$. Moreover, if $\Gamma_{f\circ g}$ denotes the graph of $f\circ g$, then we have
\[
    \Gamma_{f\circ g}=(p_{1,4})_*\mathcal C(f,g)
\]
as cycles.
\end{lemma}

\begin{proof}
The analytic set $|\cC(f,g)|$ is in fact the fibre product
\[
    |\Gamma_g|\times_X|\Gamma_f|=\left\{\bigl((x_1,x_2),(x_2,x_4)\bigr): (x_1,x_2)\in|\Gamma_g|,\ (x_2,x_4)\in|\Gamma_f|\right\},
\]
which is taken with respect to the projection maps
\[
    \pi_2:|\Gamma_g|\longrightarrow X\quad\text{and}\quad\pi_1:|\Gamma_f|\longrightarrow X.
\]

Since $f$ is holomorphic, $\pi_1$ here is finite. By \Cref{base-change} (1), its base change
\[
    |\Gamma_g|\times_X|\Gamma_f|\to |\Gamma_g|
\]
is finite as well. Consequently, every irreducible component of $|\cC(f,g)|$ has dimension at most $2$. On the other hand, $|\Gamma_g|\times|\Gamma_f|$ has pure dimension $4$, and $D$ is of dimension 6. Hence every irreducible component of $|\cC(f,g)|$ has dimension at least $2$. Thus $|\cC(f,g)|$ has pure dimension $2$, so the intersection is proper and $\cC(f,g)$ is well-defined. The identity $\Gamma_{f\circ g}=(p_{1,4})_*\mathcal C(f,g)$ follows directly from the definition of the composition of correspondences.
\end{proof}

The \textit{dynamical degree of order $q$} of $f$, denoted by $d_q$, is the spectral radius of $f^*$ on $H^{q,q}(X,\R)$. In particular, $d_0$ is the number of points of $f(x)$ for any $x\in X$, counted with multiplicity, and $d_2$ is the number of points of $f^{-1}(x)$ for a generic $x\in X$, counted with multiplicity. In the following, we will assume that $f^{-1}$ is also a holomorphic correspondence, that $d_1>\max\{d_0,d_2\}$, and that the action of $f^*$ on $H^{1,1}(X,\R)$ is \emph{simple}, i.e., $d_1$ is the only eigenvalue of maximal modulus and it is simple.

For algebraically stable meromorphic maps on surfaces, simple action is a direct consequence of the relation $d_1>d_2$. This is proven in \cite{DF01AJM}, but the proof does not work in our setting since they need the so-called push-pull formula which is not true for correspondences.

\subsection{Green currents and hyperplane sections}

Recall that correspondences under the above assumptions have an associated \emph{main Green current} $T^+$, to which smooth currents of the corresponding bidegree equidistribute under the pullback action. We normalize $T^+$ so that $\langle T^+,\omega\rangle=1$. Green currents of correspondences and equidistribution results have been studied in \cite{DNV18AM,AV24AMP,LV26}.

The following result, whose proof is adapted from \cite[Corollary 2.11]{DDG1Indiana}, gives an approximation of Green currents of correspondences on a projective surface by pullbacks of hyperplane sections. Since the associated currents are singular, the proof is more delicate than for smooth currents. In \cite[Theorem 2.1]{LV26}, the authors prove a similar result which is valid in any dimension but requires an additional assumption on the graph.

\begin{proposition}[Genericity condition (i)] \label{Gc1}
    Consider the Lebesgue measure on the dual space $\check\Pb^N$ of $\Pb^N$. For Lebesgue almost every $H$, we have $d_1^{-n}(f^n)^*[H_X]\to cT^+$, for some $c>0$ which only depends on the embedding $i:X\into\Pb^N$.
\end{proposition}

\begin{proof}
    Let $\FuStu$ be the Fubini-Study form on $\Pb^N$. By Crofton's formula \cite[Corollary 7.11]{D12Book} we have
    \begin{equation} \label{crofton}
        \FuStu=\int_{\check{\Pb}^N}[H]\d \nu(H),
    \end{equation}
    where $\nu$ denotes the Lebesgue measure on $\check{\Pb}^N$. Then for each hyperplane $H$ in $\Pb^N$, we have
    \begin{equation} \label{psi_H}
        [H]-\FuStu=\ddc\psi_H
    \end{equation}
    where $\psi_H(p)\leq 0$. For generic $H$, we can intersect both sides of \eqref{psi_H} with $X$ to obtain $[H_X]-\FuStu|_X=\ddc (\psi_H|_X)$. By \cite[Proposition 5.9]{DNV18AM}, we have $d_1^{-n}(f^n)^*(\FuStu|_X)\to cT^+$ for some $c>0$ depending only on $i$. It remains to prove that $d_1^{-n}(f^n)^*(\psi_H)\to 0$ in $L^1(X)$ for Lebesgue almost every $H$. 
    
    Notice that $(f^n)^*(\psi_H)(p)=\sum_{q\in f^n(p)} \psi_H(q)$, where the sum is counted with multiplicities. \Cref{crofton,psi_H} imply that $\int \psi_H(p)\d\nu(H)$ is pluriharmonic in $p$, hence constant since $\Pb^N$ is compact. We denote this constant by $C\leq 0$. Denote by $\d\vol_X$ the volume form on $X$ such that $\int_X\d\vol_X=1$. Therefore, by Fubini's theorem, we have
    \begin{align*}
        \int_{\check{\Pb}^N}\left\|d_1^{-n}(f^n)^*(\psi_H)\right\|_{L^1(X)}\d\nu(H)
        &=-\int_{\Check{\Pb}^N}\int_Xd_1^{-n}(f^n)^*(\psi_H)(p)\d\vol_X(p)\d\nu(H)\\
        &=-\int_X\int_{\check{\Pb}^N} d_1^{-n}(f^n)^*(\psi_H)(p)\d\nu(H) \d\vol_X(p)\\
        &=-\int_X\sum_{q\in f^n(p)}\left(\int_{\check{\Pb}^N} d_1^{-n}\psi_H(q)\d\nu(H) \right)\d\vol_X(p)\\
        &=-Cd_1^{-n}d_0^n.
    \end{align*}
    Since $d_1>d_0$, we have 
    \[
        \int_{\check{\Pb}^N}\sum_n\left\|d_1^{-n}(f^n)^*(\psi_H)\right\|_{L^1(X)}\d\nu(H)<\infty.
    \]
    Again by Fubini's theorem, for $\nu$-almost every $H$, we have $\sum_n\left\|d_1^{-n}(f^n)^*(\psi_H)\right\|_{L^1(X)}<\infty$. For such $H$, this implies $\left\|d_1^{-n}(f^n)^*(\psi_H)\right\|_{L^1(X)}\to 0$, as desired.
\end{proof}

\subsection{Woven currents} \label{woven_and_criterion} 
Woven currents were introduced in \cite{D05JGA}. We recall here the definition for $(1,1)$-currents. See also \cite{DS16MAnn} for a precise definition in higher dimension and for any bidegree.

Let $\mathcal C(X)$ denote the space of irreducible complex curves $\Gamma\subset X$, not necessarily closed in $X$, with finite volume. Let $[\Gamma]$ denote the current of integration on $\Gamma$. Let $\mu$ be a positive measure on $C(X)$. We say that $\mu$ is a \emph{web} on $X$ if $\int_{\mathcal C(X)}\vol(\Gamma)\d\mu(\Gamma)$ is finite.
    
\begin{definition}
    Let $T$ be a positive $(1,1)$-current. We say that $T$ is \emph{woven} if it can be written as
    \[
        T:=\int_{\mathcal C(X)}[\Gamma]\d\mu(\Gamma),
    \]
    where $\mu$ is a web on $X$.
\end{definition}

To check if a current is woven we will use the following form of \cite[Proposition 5.8]{Dujardin09}, see also \cite[Theorem 1]{Dujardin03MAnn} and \cite[Theorem 3.6]{DDG1Indiana}. A related local laminarity criterion for smooth analytic curves was established in \cite{dThe04ENS}. This criterion is purely geometrical and does not involve the dynamics.

We do not assume that the curves $C_n$ are irreducible, and we define the genus of a possibly reducible algebraic curve as the sum of the genera of its components. Suppose $C_n=\sum_j C_{n,j}$ where $C_{n,j}$ are the irreducible components, repeated according to their multiplicities. If ${C}_{n,j}^\nu$ is a normalization of $C_{n,j}$, then we have 
\[
    g(C_n):=\sum_j g({C}_{n,j}^\nu),
\]
which is independent of the normalization we choose.

\begin{theorem} \label{criterion}
    Let $C_n$ be a sequence of singular algebraic curves on a complex projective surface $X$, with $\deg (C_n)\to\infty$. Assume that $g(C_n)=O(\deg (C_n))$. Then the cluster values of the sequence $[C_n]/\deg(C_n)$ are woven currents. 
\end{theorem}

For a correspondence $f$ on $X$, we have the following.

\begin{corollary} \label{genus_to_woven}
    Consider the Lebesgue measure on the dual space $\check\Pb^N$ of $\Pb^N$. Then the following is true for Lebesgue almost every $H$. Set $C_n:=(f^n)^*(H_X)$. If we have $g(C_n)=O(d_1^n)$, then $T^+$ is woven.
\end{corollary}

 \begin{proof}

    By Chow's theorem, $C_n$ is an algebraic curve on $X$. Its degree with respect to the embedding $i:X\into\Pb^N$ can be computed by the intersection $\lp [C_n],\omega\rp$. By \Cref{Gc1}, we deduce that for Lebesgue almost every $H$ we have
    \[
        \lim_{n\to\infty}d_1^{-n}\deg C_n=\lp cT^+,\omega\rp=c.
    \]
    This implies
    \[
        [C_n]/\deg(C_n)\to T^+.
    \]
    Since $d_1^n=O(\deg(C_n))$, it suffices to apply \Cref{criterion} to conclude the proof.
 \end{proof}

\section{Smooth models for iterated graphs} \label{iterated_graphs}

\subsection{Construction of the smooth models}\label{subsec:construct-model}
Fix $n\geq 2$. Let $\Gamma_1$ and $\Gamma_{n-1}$ be the graphs of $f$ and $f^{n-1}$, respectively. They are effective 2-cycles in $X\times X$. In the notation of
\Cref{lem:composition-cycle}, set
    \[
        \mathcal I_n:=\mathcal C(f^{n-1},f)=D\cdot(\Gamma_1\times\Gamma_{n-1}).
    \]
Since $f^{n-1}$ is holomorphic, \Cref{lem:composition-cycle} shows that $\mathcal I_n$ is well-defined and has pure dimension $2$. Moreover,
\begin{equation}\label{Gamma-I}
    \Gamma_n=(p_{1,4})_*\mathcal I_n,
\end{equation}
where $\Gamma_n$ is the graph of $f^n$.

The following is the main result of this section.

\begin{proposition}\label{prop:resolve-graph}
    For every $n\geq 1$, there exists a triple $(\tildeGamma_n,\Pi_{1,n},\Pi_{2,n})$ such that the following are true:
    \begin{enumerate}
        \item $\tildeGamma_n$ is a smooth, possibly disconnected, compact complex surface and $\Pi_{1,n}$ and $\Pi_{2,n}$ are holomorphic maps from $\tildeGamma_n$ to $X$;
        \item both $\Pi_{1,n}$ and $\Pi_{2,n}$ are dominant and generically finite on every irreducible component of $\widetilde\Gamma_n$;
        \item  we have $(\nu_n)_*\widetilde{\Gamma}_n=\Gamma_n$, where we define $\nu_n:=(\Pi_{1,n},\Pi_{2,n})$.
    \end{enumerate}
\end{proposition}

The rest of the section is devoted to the proof of \Cref{prop:resolve-graph}. We will construct $\widetilde{\Gamma}_n$ by induction. 

\medskip
\noindent \textbf{The case of $n=1$.} Suppose $\Gamma_1=\sum_j \gamma_j$ where the $\gamma_j$'s are irreducible analytic subsets of $X\times X$, repeated according to their multiplicities. For each $j$, let $g_j:\widetilde\gamma_j\longrightarrow \gamma_j$ be a blow-up desingularization. We define
\[
    \widetilde\Gamma_1:=\bigsqcup_j \widetilde\gamma_j .
\]
Thus $\widetilde\Gamma_1$ is a smooth, possibly disconnected, complex surface. Let $\nu_1:\widetilde\Gamma_1\longrightarrow X\times X$ be the map whose restriction to $\widetilde\gamma_j$ is given by $g_j$, composed with the inclusion $\gamma_j\hookrightarrow X\times X$. We define
\[
    \Pi_{1,1}:=\pi_1\circ\nu_1\qquad\text{and}\qquad\Pi_{2,1}:=\pi_2\circ\nu_1.
\]
Then we have $\nu_1=(\Pi_{1,1},\Pi_{2,1})$ and $(\nu_1)_*\widetilde{\Gamma}_1=\Gamma_1$. By definition of a holomorphic correspondence, both $\Pi_{1,1}$ and $\Pi_{2,1}$ are dominant and generically finite on every $\widetilde{\gamma}_j$.

\medskip

\noindent \textbf{Construction of $(\tildeGamma_n,\Pi_{1,n},\Pi_{2,n})$.} Suppose that we have constructed the triple $(\widetilde{\Gamma}_{n-1},\Pi_{1,n-1},\Pi_{2,n-1})$ which satisfies the results of \Cref{prop:resolve-graph}. Set $M_n:=\widetilde{\Gamma}_1\times\widetilde{\Gamma}_{n-1}$, and let 
\[
    P_1: M_n\to \widetilde{\Gamma}_1, \qquad P_2: M_n\to \widetilde{\Gamma}_{n-1}
\]
be the two projection maps. Consider the fibre product 
\[
    |V_n|:=\tildeGamma_1\times_X\tildeGamma_{n-1}=\left\{(p,q)\in M_n:\Pi_{2,1}(p)=\Pi_{1,n-1}(q)\right\}.
\]
Let $\Delta\subset X\times X$ be the diagonal. Then it is easy to see that $|V_n|=(\Pi_{2,1}\times \Pi_{1,n-1})^{-1}(\Delta)$. On $|V_n|$, we define the map
\[
    \kappa_n:=\Pi_{2,1}\circ P_1=\Pi_{1,n-1}\circ P_2: |V_n|\to X.
\] 
Then, we have the following commutative diagram:
\[
\begin{tikzcd}[
    arrow style=tikz,
    >={Stealth[length=5pt,width=3pt]}]
{|V_n|}
  \arrow[r,"P_2"]
  \arrow[d,"P_1"']
  \arrow[dr,"\kappa_n" description, pos=0.45]
&
\widetilde{\Gamma}_{n-1}
  \arrow[d,"\Pi_{1,n-1}"]
\\
\widetilde{\Gamma}_1
  \arrow[r,"\Pi_{2,1}"']
&
X
\end{tikzcd}
\]

\begin{lemma}\label{V_n-pure-dimension 2}
    The analytic set $|V_n|$ is of pure dimension 2.
\end{lemma}

\begin{proof}
    Let $\pi:Y\longrightarrow X$ be the restriction of either $\Pi_{2,1}$ or $\Pi_{1,n-1}$ to one irreducible component of its source. By \Cref{prop:resolve-graph}, $\pi$ is proper, dominant, and generically finite, and both $Y$ and $X$ are irreducible surfaces. Set $S_\pi:=\{x\in X:\dim \pi^{-1}(x)>0\}$. We claim that $S_\pi$ is finite. Indeed, suppose that $S_\pi$ contained an irreducible curve $C$. Then $\dim \pi^{-1}(C)\ge \dim C+1=2$. Since $Y$ is irreducible and $\dim Y=2$, this would imply $\pi^{-1}(C)=Y$, hence $\pi(Y)\subseteq C$, contradicting the dominance of $\pi$. Therefore $S_\pi$ is zero-dimensional, hence it is finite. Taking the union of the sets $S_\pi$ over the finitely many irreducible components of the sources of $\Pi_{2,1}$ and $\Pi_{1,n-1}$, we obtain a finite set $S\subseteq X$ such that outside $S$ both maps have finite fibres. A similar argument gives that, over points of $S$, every fibre has dimension at most one. Now $\kappa_n^{-1}(x)=\Pi_{2,1}^{-1}(x)\times \Pi_{1,n-1}^{-1}(x)$. Hence the fibres of $\kappa_n$ are finite for $x\notin S$, and have dimension at most two for $x\in S$. Therefore $\dim |V_n|\le 2$.

    On the other hand, $|V_n|$ is the inverse image of the diagonal $\Delta\subset X\times X$, which has codimension two, under a holomorphic map from the smooth fourfold $M_n$. Thus every irreducible component of $|V_n|$ has dimension at least two. Consequently, $|V_n|$ is of pure dimension two.
\end{proof}

We can define the associated cycle on $M_n$:
\[
    V_n:=(\Pi_{2,1}\times \Pi_{1,n-1})^*\Delta.
\]
By \Cref{V_n-pure-dimension 2}, $V_n$ is an effective 2-cycle. Suppose $V_n=\sum_j W_j$ where each $W_j$ is an irreducible analytic set of pure dimension 2 in $M_n$ and they are repeated according to multiplicities. Define
\begin{equation}\label{eq: V_n^dom}
    V_n^{\mathrm{dom}}:=\sum_{\kappa_n|_{W_j} \text{ is dominant}}W_j.
\end{equation}
Then we have $V_n^{\mathrm{dom}}\leq V_n$. Finally, set
\[
    \qquad \rho:=(\nu_1\times\nu_{n-1}):M_n\to X^4.
\]

\begin{lemma}\label{lemma:V_n}
    We have $\rho_*V_n=\cI_n$. Moreover, if we denote $q:=p_{1,4}\circ \rho$, then $q_*V_n^{\mathrm{dom}}=q_*V_n=\Gamma_n$.
\end{lemma}

\begin{proof}
    Denote by $p_{2,3}:X^4 \to X\times X$ the projection which maps $(x_1,x_2,x_3,x_4)$ to $(x_2,x_3)$. Then it is easy to verify that $(\Pi_{2,1}\times\Pi_{1,n-1})=p_{2,3}\circ \rho$. Note that $D=p_{2,3}^*\Delta$. Hence $V_n=(p_{2,3}\circ \rho)^*\Delta=\rho^*D$. By the projection formula \eqref{eq:cycle-projection}, we have
    \[
        \rho_*V_n=\rho_*(\rho^*D)=D\cdot \rho_*M_n=D\cdot(\Gamma_1\times\Gamma_{n-1})=\cI_n.
    \]
    By (\ref{Gamma-I}), we have
    \[
        q_*V_n=(p_{1,4})_*(\rho_*V_n)=(p_{1,4})_*\cI_n=\Gamma_n.
    \]
    Let $W_j$ be an irreducible component of $V_n$ on which $\kappa_n$ is not dominant. Note that
    \[
        |\cI_n|=\{(x_1,x_2,x_2,x_4): (x_1,x_2)\in|\Gamma_1|, (x_2,x_4)\in|\Gamma_{n-1}|\}.
    \]
    Consider the projection map $\tau:|\cI_n|\to X$ which maps $(x_1,x_2,x_2,x_4)$ to $x_2$. Since $f^{-1}$ and $f^{n-1}$ are holomorphic, $\tau$ is a finite map, and we have $\kappa_n|_{W_j}=\tau\circ \rho|_{W_j}$. Then if $\kappa_n|_{W_j}$ is not dominant, so does $\rho|_{W_j}$. Since $q=p_{1,4}\circ \rho$, this implies $q_*W_j=0$. Therefore, we deduce that $q_*(V_n-V_n^{\mathrm{dom}})=0$.
\end{proof}

For each $W_j$ in (\ref{eq: V_n^dom}), consider its resolution of singularity $s_{n,j}:\widetilde{W}_j\to W_j$, which is an isomorphism over $W_j^{\mathrm{reg}}$. We define $\widetilde{\Gamma}_n$ to be the disjoint union of the $\widetilde{W}_j$'s, i.e.,
\begin{equation}\label{eq: tildegamma-n}
    \tildeGamma_n:=\bigsqcup_{\kappa_n|_{W_j} \text{ is dominant}} \widetilde{W}_j.
\end{equation}
Let
    \[
        s_n:=\bigsqcup s_{n,j}:\widetilde\Gamma_n\longrightarrow |V_n^{\mathrm{dom}}|
    \]
be the induced map such that for every $j$, $s_n|_{\widetilde{W}_j}:\widetilde{W}_j\to W_j$ is the resolution map. Then naturally we have 
    \[
        (s_n)_*\tildeGamma_n=V_n^{\mathrm{dom}}.
    \]
We define the following maps:
\[
    \eta_n:=P_1\circ s_n, \quad \beta_n:=P_2\circ s_n \qquad\text{and} \qquad\Pi_{1,n}:= \Pi_{1,1}\circ \eta_n,\quad \Pi_{2,n}:=\Pi_{2,n-1}\circ\beta_n.
\]
Finally, we put $\nu_n:=(\Pi_{1,n},\Pi_{2,n}):\widetilde\Gamma_n\longrightarrow X\times X$.

Then $\eta_n$ and $\beta_n$ are holomorphic maps from $\widetilde{\Gamma}_n$ to $\widetilde{\Gamma}_1$ and $\widetilde{\Gamma}_{n-1}$, respectively. By induction, $\Pi_{1,n}$ and $\Pi_{2,n}$ are holomorphic maps from $\tildeGamma_n$ to $X$. We obtain the following commutative diagram as in \cite{Dujardin03MAnn}:

\begin{equation}\label{large-pyramid}
\begin{tikzcd}[
    arrow style=tikz,
    >={Stealth[length=5pt,width=3pt]}
]
	&& {\widetilde\Gamma_n} && \\
	& {\widetilde\Gamma_1} && {\widetilde\Gamma_{n-1}} \\
	X && X && X
	\arrow["{\eta_n}"', from=1-3, to=2-2]
	\arrow["{\beta_n}", from=1-3, to=2-4]
	\arrow[bend right=50pt, "{\Pi_{1,n}}"', from=1-3, to=3-1]
	\arrow[bend left=50pt, "{\Pi_{2,n}}", from=1-3, to=3-5]
	\arrow["{\Pi_{1,1}}"'{pos=0.2}, from=2-2, to=3-1]
	\arrow["{\Pi_{2,1}}"{pos=0.2}, from=2-2, to=3-3]
	\arrow["{\Pi_{1,n-1}}"'{pos=0.2}, from=2-4, to=3-3]
	\arrow["{\Pi_{2,n-1}}"{pos=0.2}, from=2-4, to=3-5]
	\arrow["f"', from=3-1, to=3-3]
	\arrow["{f^{n-1}}"', from=3-3, to=3-5]
\end{tikzcd}
\end{equation}

We have constructed the triple $(\tildeGamma_n,\Pi_{1,n},\Pi_{2,n})$. It remains to prove (2) and (3) in \Cref{prop:resolve-graph}. This is done in the next two lemmas.

\begin{lemma} \label{nuenne}
     With the notations above, we have $\Gamma_n=(\nu_n)_*\widetilde\Gamma_n$.
\end{lemma}

\begin{proof}
    By the definitions of $\nu_n$, $\eta_n$, $\beta_n$, and $q$, we have
    \[
        \nu_n=(\Pi_{1,1}\circ \eta_n, \Pi_{2,n-1}\circ \beta_n)=(\Pi_{1,1}\circ P_1, \Pi_{2,n-1}\circ P_2)\circ s_n=p_{1,4}\circ \rho\circ s_n=q\circ s_n.
    \]
    Therefore, by \Cref{lemma:V_n}, we have
    \[
        (\nu_n)_*\tildeGamma_n=q_* (s_n)_*\tildeGamma_n=q_*V_n^{\mathrm{dom}}=\Gamma_n,
    \]
    as claimed.
\end{proof}

\begin{lemma}\label{1n-2n-dominant}
    Both $\Pi_{1,n}$ and $\Pi_{2,n}$ are dominant and generically finite on every $\widetilde{W}_j$ appeared in (\ref{eq: tildegamma-n}). Moreover, we have:
    \begin{equation*}
        \begin{aligned}
            d_t(\Pi_{1,n})=d_0^n &,\quad d_t(\Pi_{2,n})=d_2^n,\quad d_t(\eta_n)=d_0^{n-1}\quad\text{and}\quad d_t(\beta_n)=d_2.
        \end{aligned}
    \end{equation*}
\end{lemma}

\begin{proof}
Since $\kappa_n|_{W_j}=\Pi_{2,1}\circ P_1|_{W_j}$ is dominant, the image $P_1(W_j)$ has dimension two. Moreover, $P_1$ is proper, so $P_1(W_j)$ is a closed irreducible analytic subset of $\widetilde\Gamma_1$. It is therefore an irreducible component of $\widetilde\Gamma_1$. Since $\Pi_{1,1}$ is dominant on every such component and $s_{n,j}$ is surjective, the map $\Pi_{1,n}|_{\widetilde W_j}=\Pi_{1,1}\circ P_1\circ s_{n,j}$ is dominant. Similarly, the identity $\kappa_n|_{W_j}=\Pi_{1,n-1}\circ P_2|_{W_j}$ shows that $P_2(W_j)$ is an irreducible component of $\widetilde\Gamma_{n-1}$. By the induction hypothesis, $\Pi_{2,n-1}$ is dominant on this component. Hence $\Pi_{2,n}|_{\widetilde W_j}=\Pi_{2,n-1}\circ P_2\circ s_{n,j}$ is dominant as well. Both maps are generically finite, since their sources and targets have dimension two.

We next compute the degrees, taking all cycle multiplicities into account. By \Cref{nuenne} and the identities $\Pi_{1,n}=\pi_1\circ\nu_n$ and $\Pi_{2,n}=\pi_2\circ\nu_n$ we have \[
\begin{aligned}
    (\Pi_{1,n})_*\widetilde\Gamma_n=(\pi_1)_*\Gamma_n=d_0^nX,\qquad\text{and}\qquad(\Pi_{2,n})_*\widetilde\Gamma_n=(\pi_2)_*\Gamma_n=d_2^nX.
\end{aligned}
\]
This implies the first two equalities.

\medskip

\noindent\textbf{Claim.} We have $(P_i)_*V_n^{\mathrm{dom}}=(P_i)_*V_n$ for $i=1,2$.

\begin{proof}
Let $W_j$ be an irreducible component of $V_n$ on which $\kappa_n$ is not dominant. If $P_1(W_j)$ had dimension two, it would be an irreducible component of $\widetilde\Gamma_1$. Since $\Pi_{2,1}$ is dominant on every such component, the identity
\[
    \kappa_n|_{W_j}=\Pi_{2,1}\circ P_1|_{W_j}
\]
would imply that $\kappa_n|_{W_j}$ is dominant, giving a contradiction. Thus $\dim P_1(W_j)<2$, and the same argument gives $\dim P_2(W_j)<2$. Consequently, we have
\[
    (P_i)_*\bigl(V_n-V_n^{\mathrm{dom}}\bigr)=0,\qquad i=1,2.
\]
This proves the claim.
\end{proof}
Let $p$ be a generic point of an irreducible component of $\widetilde\Gamma_1$. By the definition of the fibre-product cycle $V_n$, the fibre of $P_1$ over $p$, with its natural intersection multiplicities, is identified with $\Pi_{1,n-1}^{-1}\bigl(\Pi_{2,1}(p)\bigr)$, which has total multiplicity $d_0(f^{n-1})=d_0^{n-1}$ for generic $p$. We deduce that $(P_1)_*V_n=d_0^{n-1}\widetilde\Gamma_1$. Similarly, we can prove that $ (P_2)_*V_n=d_2\widetilde\Gamma_{n-1}$. Using the equalities
\[
    (s_n)_*\widetilde\Gamma_n=V_n^{\mathrm{dom}},\qquad\eta_n=P_1\circ s_n\qquad\text{and}\qquad\beta_n=P_2\circ s_n,
\]
we deduce that
\[
\begin{aligned}
    (\eta_n)_*\widetilde\Gamma_n=(P_1)_*V_n^{\mathrm{dom}}=d_0^{n-1}\widetilde\Gamma_1,\qquad\text{and}\qquad(\beta_n)_*\widetilde\Gamma_n=(P_2)_*V_n^{\mathrm{dom}}=d_2\widetilde\Gamma_{n-1},
\end{aligned}
\]
from which we obtain the topological degrees of $\eta_n$ and $\beta_n$. This concludes the proof.
\end{proof}

\subsection{The generic hyperplane condition}

For any hyperplane $H$, we define two effective cycles
\[
    \widehat{C}_n:=\Pi_{2,n}^*(H_X), \quad C_n:=(\Pi_{1,n})_*(\widehat{C}_n),
\]
where all pushforwards are taken in the sense of cycles and thus include multiplicities. By \Cref{nuenne} and \eqref{eq:cycle-projection}, we have
\[
    \begin{aligned}
    (f^n)^*H_X=(\pi_1)_*\bigl(\pi_2^*H_X\cdot\Gamma_n\bigr)=(\pi_1)_*(\nu_n)_*\bigl(\nu_n^*(\pi_2^*H_X)\bigr)=(\Pi_{1,n})_*\Pi_{2,n}^*H_X=C_n.
    \end{aligned}
\]

In order to prove our main result, we need to choose a hyperplane $H$ satisfying some genericity conditions.

For a holomorphic map $h:X_1\to X_2$ between two smooth complex surfaces, we denote by $\cE(h)$ its exceptional locus, which is the union of irreducible curves contracted by $h$. We claim that 
\[
    \cE(\Pi_{1,n})=\cE(\Pi_{2,n})\text{ for every }n\geq 1.
\]
In fact, take an irreducible curve $C$ contracted by $\Pi_{1,n}$. Suppose $a=\Pi_{1,n}(C)\in X$. Then $\nu_n(C)\subset \pi_1^{-1}(a)\cap|\Gamma_n|$. Since $f^n$ is holomorphic, $\nu_n(C)$ must be a finite set. It is a single point because $C$ is irreducible. This implies $\Pi_{2,n}(C)$ is also a single point. Thus, $C\subset \cE(\Pi_{2,n})$ and $\cE(\Pi_{1,n})\subset\cE(\Pi_{2,n})$. The other side of the inequality can be deduced similarly.

In the following, for every $n\geq 1$, we set
\[
    \cE_n:=\cE(\Pi_{1,n})=\cE(\Pi_{2,n}).
\]

\begin{lemma}[Genericity condition (ii)] \label{Gc2}
    For Lebesgue almost every hyperplane $H\subseteq \Pb^N$, the following are true:
    \begin{enumerate}
        \item for every $n$ we have ${H_X\cap \Pi_{2,n}\big(\cE_n\big)=\emptyset}$;
        \item $\widehat{C}_n$ is a smooth, possibly disconnected, curve for every $n$.
    \end{enumerate}
\end{lemma}

\begin{proof}
    For a fixed $n$, since $\Pi_{2,n}\big(\cE_n\big)$ is finite, condition (1) is true for a generic set of hyperplanes $H$. By Bertini's theorem (see for instance \cite[Chapter 1, Section 1]{GH94PoAG}), the same is true for condition (2). For both conditions to be true for every $n$, we need to take a countable intersection of generic sets, which still gives a set of hyperplanes of full Lebesgue measure. This concludes the proof.
\end{proof}

We can now reduce ourselves to estimating the genus of the normalization. Let $\sum_j \widehat{C}_{n,j}$ be the decomposition of $\widehat{C}_n$ into irreducible components.

\begin{corollary} \label{big-o-goal_statement}
    Let $H$ be a hyperplane satisfying the conditions of \Cref{Gc1} and \Cref{Gc2}. If
    \begin{equation} \label{big-o-goal}
        G_n:=\sum_j g(\widehat{C}_{n,j})=O(d_1^n),
    \end{equation}
    then $T^+$ is woven.
\end{corollary}

\begin{proof}
    First, we notice that $\Pi_{1,n}(\widehat{C}_{n,j})$ is not a point for every $j$. Otherwise, we will have $\widehat{C}_{n,j}\subset \cE_n$. Set $a=\Pi_{2,n}(\widehat{C}_{n,j})\in \Pi_{2,n}(\cE_n)$. We also have $a\in H_X$, which contradicts \Cref{Gc2}.
    
    Set $C_{n,j}:=\Pi_{1,n}(\widehat{C}_{n,j})$. Then $(\Pi_{1,n})_*\widehat{C}_{n,j}=m_{n,j}C_{n,j}$, where $m_{n,j}\geq 1$ is the degree of $\Pi_{1,n}|_{\widehat{C}_{n,j}}$. Let $C_{n,j}^{\nu}$ be the normalization of $C_{n,j}$. By the universal property of normalizations \cite[p.164]{GR1984}, there exists a canonical map $\phi_j:\widehat{C}_{n,j}\to C_{n,j}^{\nu}$ of degree $m_{n,j}$. By Riemann-Hurwitz formula applied to $\phi_j$, we have 
    \[
        2g(\widehat C_{n,j})-2=m_{n,j}(2g(C_{n,j}^{\nu})-2)+R_{n,j},
    \]
    where $R_{n,j}\geq0$. Hence
    \[
        m_{n,j}g(C_{n,j}^{\nu})\leq g(\widehat C_{n,j})+m_{n,j}.
    \]
    Recall that $C_n=(\Pi_{1,n})_*\widehat{C}_n=\sum_j m_{n,j}C_{n,j}$, where some $C_{n,j}$'s corresponding to different $j$'s could be equal because equal images of different $\widehat{C}_{n,j}$, and thus counted with repetitions. By definition of genus, we have 
    \[
        g(C_n)=\sum_j m_{n,j} g(C_{n,j}^{\nu})\leq \sum_j g(\widehat C_{n,j})+\sum_jm_{n,j}.
    \]
    On the other hand, since every non-zero projective curve has degree at least one, we have
    \[
        \sum_jm_{n,j}\leq \sum_jm_{n,j}\deg C_{n,j}=\deg C_n=O(d_1^n),
    \]
    where the last estimate comes from the proof of \Cref{genus_to_woven}. Combined with the condition that $\sum_j g(\widehat{C}_{n,j})=O(d_1^n)$, we deduce that $g(C_n)=O(d_1^n)$. The assertion follows from \Cref{genus_to_woven}.
\end{proof}

\section{Proof of \texorpdfstring{\Cref{thm:main}}{}}\label{proof_section}

We can now prove our main theorem. Choose a reduced effective divisor $\mathcal B$ on $X$ such that
\[
    \Pi_{2,1}\bigl(\cE(\Pi_{2,1})\bigr)\cup \Pi_{2,1}\bigl(\supp R_{\Pi_{2,1}}\bigr)\subset |\mathcal B|.
\]

Set $U:=X\setminus |\cB|$. Then $U$ is a Zariski open subset of $X$.

\begin{lemma}\label{etale-pi-21}
    The restriction $\Pi_{2,1}:\Pi_{2,1}^{-1}(U)\longrightarrow U$ is a finite unramified map of degree $d_2$.
\end{lemma}
\begin{proof}
    By the choice of $U$, the Jacobian of $\Pi_{2,1}$ does not vanish on $\Pi_{2,1}^{-1}(U)$. Hence the holomorphic inverse function theorem implies that $\Pi_{2,1}:\Pi_{2,1}^{-1}(U)\to U$ is a local biholomorphism. Since $\widetilde\Gamma_1$ is compact, this map is also proper. It is therefore a finite holomorphic covering of degree $d_t(\Pi_{2,1})=d_2$.
\end{proof}

Recall that $|V_n|$ is defined to be the fibre product $\tildeGamma_1\times_X\tildeGamma_{n-1}$. Consequently, we can define
\[
|V_n|_U:=|V_n|\cap\kappa_n^{-1}(U)=\Pi_{2,1}^{-1}(U)\times_U \Pi_{1,n-1}^{-1}(U).
\]

\begin{lemma}\label{lm:base change}
    The restriction
    \[
        \beta_n:\beta_n^{-1}\bigl(\Pi_{1,n-1}^{-1}(U)\bigr)\longrightarrow\Pi_{1,n-1}^{-1}(U)
    \]
    is a finite unramified map of degree $d_2$.
\end{lemma}
    
\begin{proof}
Consider the second projection
\[
    P_2:|V_n|_U=\Pi_{2,1}^{-1}(U)\times_U\Pi_{1,n-1}^{-1}(U)\longrightarrow\Pi_{1,n-1}^{-1}(U).
\]
This map is the base change of $\Pi_{2,1}:\Pi_{2,1}^{-1}(U)\longrightarrow U$, which is a finite unramified map of degree $d_2$ by \Cref{etale-pi-21}. By \Cref{base-change}, $P_2$ is also finite unramified of degree $d_2$. In particular, $|V_n|_U$ is smooth, since $\Pi_{1,n-1}^{-1}(U)$ is smooth.

\medskip

\noindent\textbf{Claim.} $|V_n|_U\subset |V_n^{\mathrm{dom}}|$.
\begin{proof}
Let $W$ be an irreducible component of $|V_n|_U$. Since $|V_n|_U$ is smooth, $W$ is open and closed in $|V_n|_U$. The map $P_2$ is finite and open, so $P_2(W)$ is a non-empty irreducible subset that is both open and closed in $\Pi_{1,n-1}^{-1}(U)$. Hence it is an irreducible component of $\Pi_{1,n-1}^{-1}(U)$. By \Cref{1n-2n-dominant}, the map $\Pi_{1,n-1}$ is dominant on every irreducible component of $\widetilde\Gamma_{n-1}$. Therefore
\[
    \kappa_n|_W=\Pi_{1,n-1}\circ P_2|_W
\]
is dominant onto $U$. It follows that the closure of $W$ in $|V_n|$ is a component on which $\kappa_n$ is dominant. This proves the claim.
\end{proof}

For each component $W_j$ occurring in $V_n^{\mathrm{dom}}$, set $W_{j,U}:=W_j\cap |V_n|_U$. Whenever $W_{j,U}\neq\varnothing$, it is an open and closed irreducible component of the smooth space $|V_n|_U$. Consequently, the restriction $P_2|_{W_{j,U}}:W_{j,U}\longrightarrow\Pi_{1,n-1}^{-1}(U)$ is finite unramified. Moreover, $W_{j,U}\subset W_j^{\mathrm{reg}}$. By the choice of the resolution $s_{n,j}$, the restriction $s_{n,j}: s_{n,j}^{-1}(W_{j,U})\longrightarrow W_{j,U}$ is an isomorphism. Hence $P_2\circ s_{n,j}:s_{n,j}^{-1}(W_{j,U})\longrightarrow\Pi_{1,n-1}^{-1}(U)$ is finite unramified.

Since $\beta_n=P_2\circ s_n$ and $\kappa_n=\Pi_{1,n-1}\circ P_2$ on $|V_n|$, we have
\[
    \beta_n^{-1}\bigl(\Pi_{1,n-1}^{-1}(U)\bigr)=\bigsqcup_j s_{n,j}^{-1}(W_{j,U}),
\]
where the index $j$ ranges over the components of $V_n^{\mathrm{dom}}$, repeated according to their multiplicities. On each term of the above disjoint union, $\beta_n$ agrees with $P_2\circ s_{n,j}$. Thus the restriction of $\beta_n$ on the above set is a finite disjoint union of finite unramified maps, and is therefore finite unramified of degree $d_2$.
\end{proof}

\begin{lemma}[Genericity condition (iii)]\label{Gc3}
    For Lebesgue almost every hyperplane $H\subseteq \Pb^N$,  if $\widehat{C}_n=\Pi_{2,n}^*(H_X)$, then no irreducible component of $|\widehat C_n|$ is contained in $\Pi_{1,n}^{-1}(|\mathcal B|)$ for every $n$.
\end{lemma}

\begin{proof}
    Fix $n$ and let $F$ be an irreducible component of $\Pi_{1,n}^{-1}(|\mathcal B|)$. If $F\subset |\widehat C_n|=\Pi_{2,n}^{-1}(H_X)$, then $\Pi_{2,n}(F)\subset H_X$. Since $\Pi_{2,n}(F)$ is either a point or a curve, the set of hyperplanes containing it is a proper algebraic subset of the dual projective space $\check{\Pb}^N$. There are only finitely many such  components $F$ for each $n$, and $n$ ranges over a countable set. Hence the desired condition still holds for Lebesgue almost every hyperplane.
\end{proof}

We now fix a hyperplane $H$ satisfying the \hyperref[Gc1]{genericity conditions (i)}, \hyperref[Gc2]{(ii)}, and \hyperref[Gc3]{(iii)}. Since $\Pi_{2,n}=\Pi_{2,n-1}\circ\beta_n$, we have $\widehat C_n=\Pi_{2,n}^*(H_X)=\beta_n^*\widehat C_{n-1}$. Hence $\beta_n$ induces a map
\[
    b_n:=\beta_n|_{|\widehat C_n|}:|\widehat C_n|\longrightarrow|\widehat C_{n-1}|.
\]
This map is finite. In fact, if an irreducible component $D_0\subset \widehat{C}_n$ is contracted by $\beta_n$, then since $\Pi_{2,n}=\Pi_{2,n-1}\circ \beta_n$, it is also contracted by $\Pi_{2,n}$. Therefore, $D_0\subset \cE_n$. On the other hand, we also have $D_0\subset \Pi_{2,n}^{-1}(H_X)$. This implies that $H_X\cap \Pi_{2,n}(\cE_n)$ is non-empty, which contradicts \Cref{Gc2}. By \Cref{lm:base change}, $b_n$ is unramified over $ \widehat C_{n-1}\cap \Pi_{1,n-1}^{-1}(U)$. Therefore, all the branch values of $b_n$ lie in $|\widehat C_{n-1}|\cap \Pi_{1,n-1}^{-1}(|\mathcal B|)$.

Let $D$ be an irreducible component of $\widehat C_{n-1}$. By our generic choice of $H$, the curve $D$ is not contained in $\Pi_{1,n-1}^{-1}(|\mathcal B|)$. Let $D_1,\ldots,D_s$ be the irreducible components of $\widehat C_n$ which are mapped onto $D$, and put $e_i:=d_t(b_n|_{D_i})\geq 1$. Then we have $\sum_{i=1}^s e_i=d_2$. 

For every $i$, the Riemann--Hurwitz formula gives
\[
    2g(D_i)-2=e_i\bigl(2g(D)-2\bigr)+r_i,
\]
where $r_i\geq0$ is the degree of the ramification divisor of $b_n|_{D_i}$. Equivalently,
\[
    g(D_i)=e_i g(D)+1-e_i+\frac{r_i}{2}\leq e_i g(D)+\frac{r_i}{2}.
\]
Summing over $i$, we obtain
\[
    \sum_{i=1}^s g(D_i)\leq d_2g(D)+\frac12\sum_{i=1}^s r_i.
\]
Set $S_D:=D\cap \Pi_{1,n-1}^{-1}(|\mathcal B|)$. Since all branch values of the maps $b_n|_{D_i}$ belong to $S_D$, we have $\sum_{i=1}^s r_i\leq(d_2-1)\#S_D$. Consequently,
\[
    \sum_{i=1}^s g(D_i)\leq d_2g(D)+\frac{d_2-1}{2}\#S_D.
\]
Summing over all irreducible components $D$ of $\widehat C_{n-1}$ gives
\[
    G_n\leq d_2G_{n-1}+\frac{d_2-1}{2}\sum_D\#S_D.
\]
Since no component $D$ is contained in $\Pi_{1,n-1}^{-1}(|\mathcal B|)$, all intersections are proper, hence
\[
    \sum_D\#S_D\leq (\Pi_{1,n-1})^*\mathcal B\cdot\widehat C_{n-1}=\mathcal B\cdot (\Pi_{1,n-1})_*\widehat C_{n-1}=\mathcal B\cdot C_{n-1}.
\]
By the proof of \Cref{genus_to_woven}, we have $\deg C_{n-1}=O(d_1^{n-1})$, and thus $\mathcal B\cdot C_{n-1}$ is also $O(d_1^{n-1})$. We conclude that $G_n\leq d_2G_{n-1}+O(d_1^{n-1})$. By induction, we obtain
\[
    G_n\leq d_2^{\,n-1}G_1+C\sum_{k=2}^{n}d_2^{\,n-k}d_1^{\,k-1}
\]
for some constant $C>0$. Since $d_2<d_1$, we have
\[
    \sum_{k=2}^{n}d_2^{\,n-k}d_1^{\,k-1}=d_1^{n-1}\sum_{\ell=0}^{n-2}\left(\frac{d_2}{d_1}\right)^\ell=O(d_1^{n-1}).
\]
Thus $G_n=O(d_1^n)$, and \Cref{thm:main} follows from \Cref{big-o-goal_statement}.

\appendix
\refstepcounter{section}

\section*{Appendix: A family of examples}

\label[appendix]{sec:example}

In this appendix, we construct a large family of holomorphic correspondences satisfying the hypotheses of \Cref{thm:main}. The key result is given by the following proposition.

\begin{proposition}\label{example}
    Let $X$ be a smooth complex projective surface, and let $H$ be a very ample divisor on $X$. Let $\pi_1 $ and $ \pi_2 $ be the two projections $ X\times X\to X$. Choose positive integers $a,b,c,d$ such that $\max\{a,b\}\ge 2 $ and $ \max\{c,d\}\ge 2$. Then for generic choices of $D_1\in |\pi_1^*(aH)+\pi_2^*(cH)|, D_2\in |\pi_1^*(bH)+\pi_2^*(dH)|$, the set
    \[
        \Gamma:=D_1\cap D_2\subset X\times X
    \]
    defines a holomorphic correspondence $f$ on $X$ such that $f^{-1}$ is also holomorphic. The dynamical degrees of $f$ are:
    \[
        d_0=cd\,H^2,\qquad d_2=ab\,H^2,\qquad d_1=(ad+bc)H^2.
    \]
    Moreover, the eigenvalue $d_1$ is the unique eigenvalue of maximum modulus and it is simple.
\end{proposition}

Therefore, if we further choose $ad+bc>\max\{ab,cd\}$, the correspondence $f$ associated with $\Gamma$ satisfies the conditions in \Cref{thm:main}. It is also easy to check that a large subset of the lattice $\N^4$ satisfies the relations $\max\{a,b\}\geq 2$, $\max\{c,d\}\geq 2$, and $ad+bc>\max\{ab,cd\}$.

\medskip

To prove \Cref{example}, we first have the following simple lemma. 

\begin{lemma}\label{lemma:codim s+1}
    Let $C$ be an integral curve on $X$ and $s$ a positive integer. The set of $B\in |sH|=\Pb(H^0(X,\cO_X(sH)))$ such that $C\subset B$ is a linear subspace of $|sH|$ of codimension at least $s+1$.
\end{lemma}

\begin{proof}

Since $H$ is very ample, we have an embedding map
\begin{equation}\label{Kodaira}
    \Phi: x\in X\mapsto [u_0(x):\ldots:u_N(x)]\in\Pb^N
\end{equation}
for some integer $N>0$ and a set of sections $u_0,\ldots,u_N\in H^0(X,\cO(H))$ which are linearly independent. Then there exist two sections, say $u_0,u_1$, such that they are not identically 0 on $C$ and $u_0/u_1$ is a non-constant meromorphic function on $C$. Otherwise, we will have that $\Phi(C)$ is a single point, which contradicts the facts that $C$ is a curve and $\Phi$ is an embedding. Consider the sections $u_0^s,u_0^{s-1}u_1,\ldots, u_1^s$ in $H^0(X,\cO(sH))$. Their restriction to $C$ are linearly independent. In fact, suppose $\sum_{i=0}^s a_i\cdot u_0^iu_1^{s-i}=0$ on $C$ for some $a_i\in \C$. Then when $u_1\neq 0$, we have $\sum a_i (u_0/u_1)^i=0$ on $C$. Since $u_0/u_1$ is not a constant, we must have $a_i=0$ for all $i$.

Let $\tau:C\into X$ be the inclusion map. It induces the restriction map $\tau^*:H^0(X,\cO_X(sH))\to H^0(C,\cO(sH)|_C)$. By the above discussion, we deduce that
\[
    \codim\ker\tau^*=\dim \image\tau^*\geq s+1.
\]
Suppose $B=\{u_B=0\}$ for some $u_B\in H^0(X,\cO(sH))$. Then $C\subset B$ if and only if $u_B=0$ on $C$, i.e., $u_B\in\ker\tau^*$. Therefore, $C\subset B$ if and only if $B\in\Pb(\ker\tau^*)$, whose codimension in $\Pb(H^0(X,\cO_X(sH)))$ is at least $s+1$. This concludes the proof.
\end{proof}

\begin{proof}[Proof of \Cref{example}] We complete the proof in three steps. We start with the following simple observation.

\medskip
\noindent\textbf{\textit{Step 1.}} \textit{If $r,s$ are positive integers and $\max\{r,s\}\ge 2$, then the locus of pairs $(A,B)\in |rH|\times |sH|$ such that $A$ and $B$ have a common component has codimension at least $3$.}
\medskip

Without loss of generality, we may assume $s\geq 2$. For any $A\in |rH|$, let $C$ be an irreducible component of $A$. By \Cref{lemma:codim s+1}, the set of $B$ containing $C$ is of codimension at least 3 in $|sH|$. Since $A$ is a finite sum of these $C$'s, the set of $B$ having a common component with $A$ is of codimension at least 3. Since this is true for every $A\in |rH|$, this concludes Step 1.

\medskip
\noindent\textbf{\textit{Step 2.}} \textit{$f$ and $f^{-1}$ are holomorphic for generic $(D_1, D_2)$.}
\medskip

We will prove the assertion for $f$. The one for $f^{-1}$ follows from the same argument with $\pi_1$ and $\pi_2$ inverted. Suppose $D_1=\{u_1=0\}$ and $D_2=\{u_2=0\}$, for some $u_1\in H^0(X\times X,\cO_{X\times X}(\pi_1^*(aH)+\pi_2^*(cH)))$ and $u_2\in H^0(X\times X,\cO_{X\times X}(\pi_1^*(bH)+\pi_2^*(dH)))$. For any $p\in X$, the inclusion map $y\in X\mapsto (p,y)\in X\times X$ induces a linear map 
\[
    i_p^*:H^0(X\times X,\cO_{X\times X}(\pi_1^*(aH)+\pi_2^*(cH)))\to H^0(X, \cO_X(cH)),
\]
which is exactly the restriction to the fibre $\pi_1^{-1}(p)$. Since $aH$ is globally generated, we can choose $s\in H^0(X,aH)$ with $s(p)\neq0$. Given any $r\in H^0(X,cH)$, consider the section $u(x,y):=\frac{s(x)}{s(p)}\,r(y)$. Then $i_p^*u=r$. Hence $i_p^*$ is surjective.

Let $D_{1,p}:=\{i_p^*(u_1)=0\}$. When $i_p^*(u_1)$ is not identically zero, it defines an element of $|cH|$ which is exactly the fibre of $D_1$ on $\pi_1^{-1}(p)$. Similarly, we can consider the map
\[
    j_p^*: H^0(X\times X,\cO_{X\times X}(\pi_1^*(bH)+\pi_2^*(dH)))\to H^0(X, \cO_X(dH)),
\]
Let $D_{2,p}:=\{j_p^*(u_2)=0\}$. When $j_p^*(u_2)$ is not identically zero, it defines an element of $|dH|$ which is exactly the fibre of $D_2$ on $\pi_1^{-1}(p)$.  

Consider the incidence variety
\begin{align*}
    Z=\{(p,D_1,D_2)\in &X\times |\pi_1^*(aH)+\pi_2^*(cH)|\times |\pi_1^*(bH)+\pi_2^*(dH)|: \\
    &i_p^*(u_1)\equiv 0; \text{ or } j_p^*(u_2)\equiv 0; \\
    &\text{ or $D_{1,p}$ and $D_{2,p}$ are divisors on $X$ with a common component} \}.
\end{align*}

Since $H$ is very ample, we have the embedding \eqref{Kodaira}. Since $X$ is a surface, we have $N\geq 2$ and $\dim H^0(X,\cO_X(H))\geq 3$. This implies $\codim \ker i_p^*=\dim \image i_p^*\geq 3$. Similarly, we have $\codim\ker j_p^*\geq 3$. When $u_1\notin \ker i_p^*$ and $u_2\notin \ker j_p^*$, we have $D_{1,p}\in |cH|$ and $D_{2,p}\in |dH|$. By Step 1, outside a set of codimension at least 3, $D_{1,p}$ and $D_{2,p}$ have no common component. Altogether, we deduce that $Z$ has codimension at least 3. Since $X$ is a surface, the projection of $Z$ onto $|\pi_1^*(aH)+\pi_2^*(cH)|\times |\pi_1^*(bH)+\pi_2^*(dH)|$ has codimension at least 1. Hence, for generic choice of $(D_1,D_2)$, $\pi_1|_{\Gamma}^{-1}(p)=D_{1,p}\cap D_{2,p}$ is a finite set of points for every $p\in X$. Finally, since $\Gamma$ is locally defined by two equations, every irreducible component $Y$ of $\Gamma$ has dimension at least 2. Therefore, $\pi_i(Y)$ must be the whole $X$ for $i=1,2$, otherwise $\pi_i|_Y$ will have a fibre of positive dimension. This implies that $f$ is holomorphic and concludes Step 2.

\medskip
\noindent\textbf{\textit{Step 3.}} \textit{Computation of the dynamical degrees.}
\medskip

Let $h=c_1(\mathcal O_X(H))$, and write $h_1=\pi_1^*h $ and $ h_2=\pi_2^*h$. Then we have $\{D_1\}=ah_1+ch_2 $ and $ \{D_2\}=bh_1+dh_2$. Hence,
\[
    \{\Gamma\}=\{D_1\}\cdot\{D_2\}=ab\,h_1^2+(ad+bc)h_1h_2+cd\,h_2^2.
\]

For a generic point $x\in X$, the fibre $\pi_1|_{\Gamma}^{-1}(x)$ is the intersection of two divisors in $|cH|$ and $|dH|$. Thus $d_0=\deg(\pi_1|_\Gamma)=cd\,H^2$. Similarly, we have $d_2=\deg(\pi_2|_\Gamma)=ab\,H^2$. For $\alpha\in H^{1,1}(X,\mathbb R)$, the pullback induced by the correspondence is
\[
    f^*\alpha=(\pi_1)_*\bigl(\{\Gamma\}\smile \pi_2^*\alpha\bigr).
\]
Notice that in the expression of $\{\Gamma\}$, only the middle term $(ad+bc)h_1h_2$ contributes in the above formula. We obtain
\begin{align*}
    f^*\alpha&=(ad+bc)(\pi_1)_*(h_1h_2\smile \pi_2^*\alpha)=(ad+bc)(\pi_1)_*\big(\pi_1^*(h)\smile \pi_2^*(h\smile\alpha)\big)\\
    &=(ad+bc)\left(\int_X H\wedge\alpha\right)\cdot h
    =(ad+bc)(H\cdot\alpha)h.
\end{align*}
In particular, $f^*h=(ad+bc)H^2h$, and the image of $f^*$ is of dimension one. Therefore the eigenvalues of $f^*$ on $H^{1,1}(X,\mathbb R)$ are $(ad+bc)H^2,0,\ldots,0$. Thus $d_1=(ad+bc)H^2$, and it is a simple eigenvalue. This concludes Step 3 and the proof of \Cref{example}.
\end{proof}

\begin{example}
    Take $X=\Pb^2$ and $H$ any line in $\Pb^2$ as in \Cref{example}. For $x=(x_0,x_1,x_2)\in \C^3$ and $I=(i_0,i_1,i_2)\in \N^3$, denote by $x^I=x_0^{i_0}x_1^{i_1}x_2^{i_2}$. Let $|I|=i_0+i_1+i_2$. For positive integers $a$, $b$, $c$, and $d$, and constants $A_{I,J},B_{I,J}\in\C$, let 
\begin{align*}
    F_1(x,y)=\sum_{|I|=a,|J|=c} A_{I,J}x^Iy^J,\qquad F_2(x,y)=\sum_{|I|=b,|J|=d} B_{I,J}x^I y^J
\end{align*}
be two homogeneous polynomials in $x$ and $y$. Then for a general choice of the coefficients $A_{I,J},B_{I,J}$, $\{F_1=0\}$ and $\{F_2=0\}$ are divisors in $|\pi_1^*(aH)+\pi_2^*(cH)|$ and $|\pi_1^*(bH)+\pi_2^*(dH)|$, respectively. By \Cref{example}, when $ad+bc>\max\{ab,cd\}$, $\max\{a,b\}\geq 2$, and $\max\{c,d\}\geq 2$, the correspondence $f:[x]\mapsto [y]$ defined by 
\[
    \Gamma:=\{([x],[y])\in\Pb^2\times\Pb^2: F_1(x,y)=F_2(x,y)=0\}
\]
is holomorphic and satisfies the conditions in \Cref{thm:main}. 
\end{example}


\printbibliography[heading=bibliography,title={References}]

\end{document}